\documentclass[a4paper,12pt]{amsart}
\usepackage{amssymb,amsmath,a4wide,graphicx,stmaryrd,fullpage,setspace,microtype,textcomp,tikz,tikz-cd,accents,mathtools,abraces}
\usepackage{hyperref}
\usepackage[shortlabels]{enumitem}
\definecolor{dark-red}{rgb}{0.5,0.15,0.15}
\hypersetup{
	colorlinks   = true, 
	urlcolor     = dark-red, 
	linkcolor    = black, 
	citecolor   = dark-red 
}
\usepackage[T1]{fontenc}
\usepackage{lmodern}
\usepackage[numbers,sort&compress]{natbib}

\title{Directed paths and Moore flows}

\author[P. Gaucher]{Philippe Gaucher}

\address{Universit\'e Paris Cit\'e, CNRS, IRIF, F-75013, Paris, France}

\urladdr{\url{http://www.irif.fr/~gaucher}} 

\makeatletter
\@namedef{subjclassname@2020}{\textup{2020} Mathematics Subject Classification}
\makeatother


\swapnumbers

\makeatletter
\let\leq\@undefined
\let\geq\@undefined
\let\vec\@undefined
\makeatother
\newcommand{\leq}{\leqslant}
\newcommand{\geq}{\geqslant}
\newcommand{\vec}{\overrightarrow}

\newtheorem{thm}{Theorem}[section]
\newtheorem{prop}[thm]{Proposition}
\newtheorem{lem}[thm]{Lemma}
\newtheorem{cor}[thm]{Corollary}

\newtheorem{defnot}[thm]{Definition and notation}
\newcommand{\bdn}{\begin{defnot}}
\newcommand{\edn}{\end{defnot}}
\newcommand{\bp}{\begin{prop}}
\newcommand{\ep}{\end{prop}}
\newcommand{\bth}{\begin{thm}}
\renewcommand{\eth}{\end{thm}}
\newcommand{\bpf}{\begin{proof}}
\newcommand{\epf}{\end{proof}}
\newcommand{\bc}{\begin{cor}}
\newcommand{\ec}{\end{cor}}

\theoremstyle{definition}
\newtheorem{defn}[thm]{Definition}

\newcommand{\bd}{\begin{defn}}
\newcommand{\ed}{\end{defn}}
\newtheorem{nota}[thm]{Notation}

\newcommand{\C}{\mathcal{C}}
\newcommand{\W}{\mathcal{W}}
\newcommand{\F}{\mathcal{F}}
\newcommand{\de}{\partial}
\newcommand{\p}{\times}
\newcommand{\PP}{\mathbb{P}}
\newcommand{\Top}{{\mathbf{Top}}}
\newcommand{\iso}{\cong}
\newcommand{\ot}{\otimes}
\newcommand{\ptop}[1]{{\brm{{#1}dTop}}}
\newcommand{\moore}{{\mathbb{M}}}
\newcommand{\tto}{\simeq^+}
\newcommand{\brm}[1]{{\rm{\mathbf{#1}}}}

\newcommand{\dtopM}{{\mathcal{M}\brm{Flow}}}
\newcommand{\set}{{\brm{Set}}}
\newcommand{\ttop}{{\brm{TOP}}}
\newcommand{\Topdgr}{[\mathcal{M}^{op},\Top]}
\newcommand{\Topdgrq}{[\mathcal{M}^{op},\Top_q]_0^{proj}}
\newcommand{\Topdgrm}{[\mathcal{M}^{op},\Top_m]_0^{proj}}
\newcommand{\Topdgrh}{[\mathcal{M}^{op},\Top_h]_0^{proj}}
\DeclareMathOperator{\id}{Id}
\newcommand{\liminj}{\varinjlim}
\DeclareMathOperator{\natcub}{\underline{nat}^\square}

\allowdisplaybreaks

\begin{document}

\begin{abstract}
	This addendum extends prior work to the non-regular setting by introducing the tame realization of a precubical set as a multipointed $d$-space. Its execution paths are precisely the nonconstant tame $d$-paths in the geometric realization of the precubical set. The associated Moore flow induces a functor from precubical sets to Moore flows, which is naturally weakly equivalent, within the $h$-model structure, to a colimit-preserving functor whose image is included in the class of m-cofibrant Moore flows. For spatial (and thus proper) precubical sets, these functors coincide.
\end{abstract}
	
\maketitle
\hypersetup{linkcolor = dark-red}

\section{Introduction}

This addendum extends the results of \cite{RegularMoore} to the broader setting of not necessarily regular $d$-paths. The key modification is replacing the reparametrization category $\mathcal{G}$, whose set of objects is the open interval $]0,+\infty[$ and morphisms $\mathcal{G}(\ell_1,\ell_2)$ are nondecreasing homeomorphisms from $[0,\ell_1]$ to $[0,\ell_2]$ \cite[Proposition~4.9]{Moore1}, with the reparametrization category $\mathcal{M}$, which retains the same objects but whose morphisms $\mathcal{M}(\ell_1,\ell_2)$ are nondecreasing surjective maps from $[0,\ell_1]$ to $[0,\ell_2]$ \cite[Proposition~4.11]{Moore1}.

This transition parallels the shift from \cite{Moore2}, which establishes a Quillen equivalence between multipointed $d$-spaces and Moore flows for the reparametrization category $\mathcal{G}$, to \cite{Moore3}, extending this Quillen equivalence to the reparametrization category $\mathcal{M}$. Since \cite{Moore3} was not available when \cite{RegularMoore} was written, this addendum could not have been written earlier. We refer to \cite{RegularMoore} for identical notions and notations. Most proofs adapt by replacing $\mathcal{G}$ with $\mathcal{M}$, leveraging the categorical properties of reparametrization categories from \cite{Moore1,Moore2,Moore3}. The only significant modification is the proof of Theorem~\ref{thm:Psi} (cf. \cite[Theorem~4.11]{RegularMoore}).

By equipping the topological $n$-cube with its set of nonconstant $d$-paths between two vertices, we introduce the tame realization $|K|^t$ of a precubical set $K$ as a multipointed $d$-space (Definition~\ref{def:tame_rea}). All execution paths of $|K|^t$ are nonconstant tame $d$-paths in the sense of \cite[Section~2.9]{MR4070250}. By the tamification theorem \cite[Theorem~6.1]{MR4070250}, this point of view is not restrictive.

For concurrency theory, the homotopy type of the space of execution paths, not the topology of the underlying state space, is critical \cite{DAT_book}. Thus, we focus on the Moore flow $\moore^{top}(|K|^t)$, obtained by retaining only the execution paths and their reparametrization data from $|K|^t$ (Definition~\ref{def:Moore-flow} and Theorem~\ref{thm:Moore_top}).

The functor $K \mapsto \moore^{top}(|K|^t)$ is the composition of a left adjoint $K \mapsto |K|^t$ from precubical sets to multipointed $d$-spaces and a right adjoint $X\mapsto\moore^{top}(X)$ from multipointed $d$-spaces to Moore flows. Though this composite fails to preserve colimits in general, our main result (Definition~\ref{def:colimit-tame-rea}, Theorem~\ref{thm:iso_reg_reg0}) replaces it --up to a natural weak equivalence in the $h$-model structure-- with a colimit-preserving functor $K\mapsto [K]^{t}$, ensuring homotopy-equivalent execution path spaces, whose image is included in the class of $m$-cofibrant Moore flows (Theorem~\ref{thm:m-cof-moore-flow}). For spatial precubical sets (Definition~\ref{def:carac_spatial}), the two functors coincide, yielding an isomorphism $\liminj \moore^{top}(|K_i|^t) \iso \moore^{top}(|\liminj_i K_i|^t)$ for colimits $\liminj_i K_i$ of spatial precubical sets with spatial colimits. Note that this statement is optimal.

All large categories in this note are locally presentable, and all model categories are at least accessible, and sometimes combinatorial.

\section{The tame realization of a precubical set}

Unless otherwise specified, symbols such as $\ell$, $\ell'$, $\ell_i$, $L$, and similar denote strictly positive real numbers. Let $\mu_\ell: [0,\ell] \to [0,1]$ be the homeomorphism defined by $\mu_\ell(t) = t/\ell$.

The notation $[0,\ell_1] \tto [0,\ell_2]$ refers to a \emph{nondecreasing surjective map} from $[0,\ell_1]$ to $[0,\ell_2]$, where $\ell_1, \ell_2 > 0$. The enriched small category $\mathcal{M}$ is defined as follows: 
\begin{itemize}
	\item The set of objects is the open interval $]0,+\infty[$.
	\item The space $\mathcal{M}(\ell_1,\ell_2)$ is the set $\{[0,\ell_1]\tto [0,\ell_2]\}$ for all $\ell_1,\ell_2>0$ equipped with the $\Delta$-kelleyfication of the compact-open topology.
	\item For every $\ell_1,\ell_2,\ell_3>0$, the composition map \[\mathcal{M}(\ell_1,\ell_2)\p \mathcal{M}(\ell_2,\ell_3) \to \mathcal{M}(\ell_1,\ell_3)\] is induced by the composition of continuous maps. 
\end{itemize}

\bd \cite[Definition~6]{Moore3} A \emph{multipointed $d$-space $X$}  is a triple $(|X|,X^0,\PP^{top}X)$ where
\begin{itemize}
	\item The pair $(|X|,X^0)$ is a multipointed space, i.e. a topological space $|X|$ and a subset $X^0\subset |X|$. The space $|X|$ is called the \emph{underlying space} of $X$ and the set $X^0$ the \emph{set of states} of $X$.
	\item The set $\PP^{top}X$ is a set of continuous maps from $[0,1]$ to $|X|$ called the \emph{execution paths}, satisfying the following axioms:
	\begin{itemize}
		\item For any execution path $\gamma$, one has $\gamma(0),\gamma(1)\in X^0$.
		\item Let $\gamma$ be an execution path of $X$. Then any composite $\gamma\phi$ with $\phi:[0,1]\tto [0,1]$ is an execution path of $X$.
		\item Let $\gamma_1$ and $\gamma_2$ be two composable execution paths of $X$; then the normalized composition $\gamma_1 *_N \gamma_2$, defined by $(\gamma_1 *_N \gamma_2)(t)=\gamma_1(2t)$ for $0\leq t\leq 1/2$ and $(\gamma_1 *_N \gamma_2)(t)=\gamma_2(2t-1)$ for $1/2\leq t\leq 1$,  is an execution path of $X$.
	\end{itemize}
\end{itemize}
A map $f:X\to Y$ of multipointed $d$-spaces is a map of multipointed spaces from $(|X|,X^0)$ to $(|Y|,Y^0)$ such that for any execution path $\gamma$ of $X$, the map $\PP^{top}f:\gamma\mapsto f. \gamma$ is an execution path of $Y$. The category of multipointed $d$-spaces is denoted by $\ptop{\mathcal{M}}$. 
\ed

The subset of execution paths from $\alpha$ to $\beta$ is the set of $\gamma\in\PP^{top} X$  such that $\gamma(0)=\alpha$ and $\gamma(1)=\beta$; it is denoted by $\PP^{top}_{\alpha,\beta} X$. The set $\PP^{top}_{\alpha,\beta} X$ is equipped with the $\Delta$-kelleyfication of the compact-open topology.

\bp \label{prop:cocubical-tame-cube} (compare with \cite[Proposition~2.22]{RegularMoore})
Let $n\geq 1$. The following data assemble into a multipointed $d$-space $|\square[n]|^{t}$ called the \emph{tame $n$-cube}:
\begin{itemize}
	\item The underlying space is the topological $n$-cube $[0,1]^n$.
	\item The set of states is $\{0,1\}^n\subset [0,1]^n$.
	\item The set of execution paths from $\underline{a}$ to $\underline{b}$ with $\underline{a}<\underline{b} \in \{0<1\}^n$ is the set of paths $[0,1]\to [0,1]^n$ from $\underline{a}$ to $\underline{b}$ which are nondecreasing with respect to each axis of coordinates.
	\item The set of execution paths from $\underline{a}$ to $\underline{b}$ with $\underline{a}\geq \underline{b}$ is empty.
\end{itemize}
Let $|\square[0]|^{t} = \{()\}$. 	Let $\delta_i^\alpha : [0,1]^{n-1} \rightarrow [0,1]^n$ be the continuous map defined for $1\leq i\leq n$ and $\alpha \in \{0,1\}$ by $\delta_i^\alpha(\epsilon_1, \dots, \epsilon_{n-1}) = (\epsilon_1,\dots, \epsilon_{i-1}, \alpha, \epsilon_i, \dots, \epsilon_{n-1})$. By convention, let $[0,1]^0=\{()\}$. The mapping $[n]\mapsto |\square[n]|^{t}$ yields a well-defined cocubical multipointed $d$-space.
\ep

\bd \label{def:tame_rea} (compare with \cite[Definition~2.23]{RegularMoore})
Let $K$ be a precubical set. The \emph{tame realization} of $K$ is the multipointed $d$-space 
\[
|K|^{t} = \liminj_{\square[n]\rightarrow K} |\square[n]|^{t}.
\]
This yields a colimit-preserving functor from precubical sets to multipointed $d$-spaces.
\ed

Since the underlying space functor from multipointed $d$-spaces to spaces is colimit-preserving, the underlying space of the multipointed $d$-space $|K|^{t}$ is the topological space 
\[
|K|_{geom} = \liminj_{\square[n]\to K} [0,1]^n.
\]
Let $\underline{x}=(x_1,\dots,x_n)$ and $\underline{x}'=(x'_1,\dots,x'_n)$ be two elements of $[0,1]^n$ with $n\geq 1$. Let
	\[
	d_1(\underline{x},\underline{x}') = \sum_{i=1}^{n} |x_i-x'_i|.
	\]
The $d_1$ metric on cubes extends to a pseudometric space $|K|_{d_1}$ on the underlying set of $|K|_{geom}$, as the category of pseudometric spaces is bicomplete \cite[Definition~4.6]{RegularMoore}. The homeomorphism $|\square[n]|_{geom} \iso |\square[n]|_{d_1}$ induces a continuous bijection $|K|_{geom} \to |K|_{d_1}$, which is a homeomorphism if and only if $K$ is locally finite in a suitable sense.

\section{$L_1$-arc length of $d$-paths of precubical sets}
\label{length_sec}

\bd \label{def:path-with-length} (compare with \cite[Definition~4.1]{RegularMoore})
Let $K$ be a precubical set. Let $(\alpha,\beta)\in K_0\p K_0$. Let $\ell>0$. Let $\vec{P}^\ell_{\alpha,\beta}(K)$ be the subspace of continuous maps from $[0,\ell]$ to $|K|_{geom}$ defined by \[\vec{P}^\ell_{\alpha,\beta}(K) = \{\gamma\mu_\ell\mid \gamma\in \PP^{top}_{\alpha,\beta}|K|^{t}\}.\] 
Its elements are called \emph{the tame $d$-paths of $K$ (of length $\ell$)} from $\alpha$ to $\beta$. Let 
\[
\vec{P}^\ell(K) = \coprod_{(\alpha,\beta)\in K_0\p K_0} \vec{P}^\ell_{\alpha,\beta}(K).
\]
\ed

A fundamental property of all $d$-paths in the geometric realization of a precubical set $K$ is the existence of a well-defined $L_1$-arc length \cite[Section~2.2.1]{MR2521708} \cite[Section~2.2]{Raussen2012}. For a $d$-path $\gamma: [0,\ell] \to [0,1]^n$, the $L_1$-arc length between $\gamma(t)$ and $\gamma(t')$ is given by the distance $d_1(\gamma(t), \gamma(t'))$. Consequently, the $L_1$-arc length of a $d$-path in $[0,1]^n$ from $(0,0,\dots,0)$ to $(1,1,\dots,1)$ is $n$. For a Moore composition of $d$-paths, the $L_1$-arc length is simply the sum of the $L_1$-arc lengths of its constituent $d$-paths.

\bp \label{prop:same-length} (compare with \cite[Proposition~4.5]{RegularMoore})
Let $K$ be a precubical set. Let $(\alpha,\beta)\in K_0\p K_0$. Two execution paths of $\PP_{\alpha,\beta}^{top}|K|^{t}$ which are in the same path-connected component have the same $L_1$-arc length. 
\ep

\bpf
Write $(\PP_{\alpha,\beta}^{top}|K|^{t})_{co}$ for the underlying set of $\PP_{\alpha,\beta}^{top}|K|^{t}$ endowed with the compact-open topology. Let $\vec{P}_{\alpha,\beta}(K)_{co}$ be the set of \emph{all} $d$-paths of $K$ from $\alpha$ to $\beta$ equipped with the compact-open topology. By \cite[Proposition~2.2]{Raussen2012}, a composite map of the form $[0,1]\to (\PP_{\alpha,\beta}^{top}|K|^{t})_{co} \subset \vec{P}_{\alpha,\beta}(K)_{co} \to \mathbb{R}$ is constant, where the right-hand map is the $L_1$-arc length. The $\Delta$-kelleyfication functor does not change the path-connected components. Hence the proof is complete. 
\epf

\bp \label{prop:cont} (compare with \cite[Proposition~4.7]{RegularMoore})
Let $K$ be a precubical set. Let $\ell >0$. The set map \[L:\vec{P}^\ell(K) \p [0,\ell] \longrightarrow [0,+\infty[\] which takes $(\gamma,t)$ to the $L_1$-arc length between $\gamma(0)$ and $\gamma(t)$ is continuous. By adjunction, we obtain a continuous map \[\widehat{L}:\vec{P}^\ell(K) \longrightarrow \ttop([0,\ell],[0,+\infty[).\] 
\ep

\bpf
Let $\vec{P}(K)_{d_1}$ be the set of all $d$-paths of $K$ equipped with the compact-open topology associated with the underlying topological space of $|K|_{d_1}$. Using \cite[Lemma~2.13]{MR2521708}, the set map $L:\vec{P}(K)_{d_1} \p [0,+\infty[ \longrightarrow [0,+\infty[$ taking a pair $(\gamma,t)$ to the $L_1$-arc length between $\gamma(0)$ and $\gamma(t)$ is continuous. Let $(\vec{P}^\ell(K))_{d_1}$ ($(\vec{P}^\ell(K))_{co}$ resp.) be the underlying set of the space $\vec{P}^\ell(K)$ equipped with the compact-open topology associated with the underlying topological space of $|K|_{d_1}$ (associated with the topological space $|K|_{geom}$ resp.). Since the identity induces a continuous map $|K|_{geom} \to |K|_{d_1}$, we obtain a continuous map 
\[L:(\vec{P}^\ell(K))_{co} \p [0,\ell] \to (\vec{P}^\ell(K))_{d_1} \p [0,\ell] \subset \vec{P}(K)_{d_1} \p [0,+\infty[\to [0,+\infty[.\]
Finally, take the image by the $\Delta$-kelleyfication functor. The latter is a right adjoint therefore it preserves binary products. Besides, $[0,\ell]$ and $[0,+\infty[$ are already $\Delta$-generated. Hence the proof is complete.
\epf

\begin{lem} \label{lem:minicalcul} (compare with \cite[Lemma~4.9]{RegularMoore})
	For all $r\in \vec{P}^\ell(K)$, for all $\phi\in \mathcal{M}(\ell',\ell)$, and for all $t\in [0,\ell']$, one has \[\widehat{L}(r\phi)(t) = \widehat{L}(r)(\phi(t)).\]
\end{lem}

\bpf
The $L_1$-arc length between $r(0)=r(\phi(0))$ and $r(\phi(t))$ for the $d$-path $r$ is equal to the $L_1$-arc length between $(r\phi)(0)$ and $(r\phi)(t)$ for the $d$-path $r\phi$.
\epf

Intuitively, the natural $d$-paths are the $d$-paths whose speed corresponds to the $L_1$-arc length.

\bd \cite[Definition~2.14]{MR2521708}
Let $K$ be a precubical set. Let $(\alpha,\beta)\in K_0\p K_0$. A $d$-path $\gamma$ of $\vec{P}^\ell_{\alpha,\beta}(K)$ is \emph{natural} if $\widehat{L}(\gamma)(t)=t$ for all $t\in [0,\ell]$. This implies that $\ell$ is an integer (greater than or equal to $1$). The subset of natural $d$-paths of length $n\geq 1$ from $\alpha$ to $\beta$ equipped with the $\Delta$-kelleyfication of the compact-open topology is denoted by $\vec{N}^n_{\alpha,\beta}(K)$. 
\ed

The following theorem is an improvement of \cite[Proposition~2.16]{MR2521708}: the homotopy equivalence is replaced by a homeomorphism thanks to the spaces $\mathcal{M}(\ell,n)$. 

\bth \label{thm:Psi} (compare with \cite[Theorem~4.11]{RegularMoore})
Let $K$ be a precubical set. Let $(\alpha,\beta)\in K_0\p K_0$. The continuous map $\Phi^{\ell}:(\phi,\gamma)\mapsto \gamma\phi$ from $\mathcal{M}(\ell,n) \p \vec{P}_{\alpha,\beta}^n(K)$ to $\vec{P}_{\alpha,\beta}^{\ell,n}(K)$ induces a homeomorphism for all $n\geq 1$, and therefore a homeomorphism 
\[
\Phi^{\ell}: \displaystyle\coprod_{n\geq 1} \mathcal{M}(\ell,n) \p \vec{N}_{\alpha,\beta}^n(K) \stackrel{\iso}\longrightarrow \vec{P}_{\alpha,\beta}^{\ell}(K).
\]
\eth

\bpf
By Proposition~\ref{prop:same-length}, the space $\vec{P}_{\alpha,\beta}^{\ell}(K)$ is the direct sum of the subspaces $\vec{P}_{\alpha,\beta}^{\ell,n}(K)\supset \vec{N}_{\alpha,\beta}^{n}(K)$ of tame $d$-paths of $L_1$-arc length $n$ for $n\geq 1$. It then suffices to prove the homeomorphism $\mathcal{M}(\ell,n) \p \vec{N}_{\alpha,\beta}^n(K) \iso \vec{P}_{\alpha,\beta}^{\ell,n}(K)$ for all $n\geq 1$.

The mapping $\natcub:\gamma\mapsto \gamma\widehat{L}(\gamma)^{-1}$ (called the \emph{naturalization}: see \cite[Definition~2.14]{MR2521708}) is a well-defined continuous map from $\vec{P}_{\alpha,\beta}^{\ell,n}(K)$ to $\vec{N}_{\alpha,\beta}^n(K)$ by \cite[Proposition~2.15]{MR2521708} for the compact-open topology. By applying the $\Delta$-kelleyfication functor which does not change the underlying set, we obtain the continuity of the naturalization map for the topologies used in this note. 

The mapping $\Psi^\ell:\gamma\mapsto (\widehat{L}(\gamma),\gamma\widehat{L}(\gamma)^{-1})$ from $\vec{P}_{\alpha,\beta}^{\ell,n}(K)$ to $\mathcal{M}(\ell,n) \p \vec{N}_{\alpha,\beta}^n(K)$ is then continuous by Proposition~\ref{prop:cont}. One has $\Phi^\ell\Psi^\ell=\id_{\vec{P}_{\alpha,\beta}^{\ell,n}(K)}$. Thus $\Phi^\ell$ is surjective. 

Assume that $\Phi^\ell(\phi_1,\gamma_1)=\Phi^\ell(\phi_2,\gamma_2)=\gamma$. Then $\gamma_1=\natcub(\gamma)=\gamma_2$. By \cite[Proposition~19]{Moore3} which is a Moore variant of \cite[Lemma~3.9]{reparam}, we deduce that $\phi_1=\phi_2$, $\natcub(\gamma)$ being regular and $|K|_{geom}$ being Hausdorff. Therefore $\Phi^\ell$, as well as $\Psi^\ell$, are bijective and the proof is complete.
\epf

\section{The \{q,h,m\}-model structures of Moore flows}

\bd \label{def:mspace} \cite[Definition~5.1]{Moore1}
The category of contravariant functors from $\mathcal{M}$ to $\Top$ is denoted by $\Topdgr_0$. The objects of $\Topdgr_0$ are called the \emph{$\mathcal{M}$-spaces}. 
\ed

\bth (\cite[Theorem~5.14]{Moore1}) \label{thm:closedsemimonoidal}
Let $D$ and $E$ be two $\mathcal{M}$-spaces. Let 
\[
D \ot E = \int^{(\ell_1,\ell_2)} \mathcal{M}(-,\ell_1+\ell_2) \p D(\ell_1) \p E(\ell_2).
\]
The pair $(\Topdgr_0,\ot)$ has the structure of a biclosed semimonoidal category.
\eth

By \cite[Theorem~6.5(ii)]{MoserLyne}, since all topological spaces are fibrant and since $(\Top,\p,\{0\})$ is a locally presentable base by \cite[Corollary~3.3]{dgrtop}, the category of $\mathcal{M}$-spaces $\Topdgr_0$ can be endowed with the projective model structure associated with one of the three model structures $\Top_q$ \cite[Theorem~2.4.19]{MR99h:55031}, $\Top_h$ \cite[Theorem~3]{vstrom3} or $\Top_m$ \cite[Theorem~2.1]{mixed-cole}. They are called the projective q-model structure (h-model structure, m-model structure resp.) and denoted by $\Topdgrq$ ($\Topdgrh$, $\Topdgrm$ resp.). All $\mathcal{M}$-spaces are fibrant for these three model structures.

\bd \cite[Definition~6.2]{Moore1} \label{def:Moore-flow}
A \emph{Moore flow} is a small semicategory enriched over the biclosed semimonoidal category $(\Topdgr_0,\ot)$ of Theorem~\ref{thm:closedsemimonoidal}. The corresponding category is denoted by $\dtopM$. 
\ed

For a Moore flow $X$, the topological space $\PP_{\alpha,\beta}X(\ell)$ is denoted by $\PP_{\alpha,\beta}^\ell X$.

\bth (\cite[Theorem~3]{Moore3}) \label{thm:Moore_top}
Let $X$ be a multipointed $d$-space. Let $\PP^\ell_{\alpha,\beta}X$ be the subspace of continuous maps from $[0,\ell]$ to $|X|$ defined by $\PP^\ell_{\alpha,\beta}X = \{\gamma\mu_\ell\mid \gamma\in \PP^{top}_{\alpha,\beta}X\}$. The following data assemble into a Moore flow $\moore^{top}(X)$:
\begin{itemize}
	\item The set of states $X^0$ of $X$
	\item For all $\alpha,\beta\in X^0$ and all real numbers $\ell>0$, $\PP_{\alpha,\beta}^{\ell}\moore^{top}(X) = \PP_{\alpha,\beta}^{\ell}X$.
	\item For all maps $[0,\ell]\tto[0,\ell']$, a map $f:[0,\ell']\to |X|$ of $\PP_{\alpha,\beta}^{\ell'}\moore^{top}(X)$ is mapped to the map $[0,\ell]\tto[0,\ell']\stackrel{f}\to |X|$ of $\PP_{\alpha,\beta}^{\ell}\moore^{top}(X)$ 
	\item For all $\alpha,\beta,\gamma\in X^0$ and all real numbers $\ell,\ell'>0$, the composition maps \[*:\PP_{\alpha,\beta}^{\ell}\moore^{top}(X) \p \PP_{\beta,\gamma}^{\ell'}\moore^{top}(X) \to \PP_{\alpha,\gamma}^{\ell+\ell'}\moore^{top}(X)\] is the Moore composition.
\end{itemize}
The mapping $\moore^{top}$ induces a functor from multipointed $d$-spaces to Moore flows which is a right adjoint. The left adjoint is explicitly described in \cite[Appendix~B]{Moore2}.   
\eth

\bth \label{thm:qhmMooreFlow}
Let $(\C,\F,\W)$ be either the projective q-model structure, or the projective h-model structure, or the projective m-model structure of $\mathcal{M}$-spaces. Then the category of Moore flows can be endowed with an accessible model structure characterized as follows: 
\begin{itemize}
	\item A map of Moore flows $f:X\to Y$ is a weak equivalence if and only if $f^0:X^0\to Y^0$ is a bijection and $\PP f:\PP_{\alpha,\beta}X\to \PP_{f(\alpha),f(\beta)}Y$ belongs to $\W$.
	\item A map of Moore flows $f:X\to Y$ is a fibration if and only if $\PP f:\PP_{\alpha,\beta}X\to \PP_{f(\alpha),f(\beta)}Y$ belongs to $\F$.
\end{itemize}
All objects are fibrant. The m-model structure of Moore flows is the mixing of their q-model structure and h-model structure in the sense of \cite[Theorem~2.1]{mixed-cole}. The q-model structure coincides with the one of \cite[Theorem~8.8]{Moore1}. These three model categories are denoted by $\dtopM_q,\dtopM_h,\dtopM_m$.
\eth

\bpf
It is mutatis mutandis the proof of \cite[Theorem~3.19]{RegularMoore}.
\epf

\section{Colimit-preserving tame realization}

\bd \label{def:colimit-tame-rea}
The cocubical Moore flow $\moore^{top}(|\square[*]|^{t})$ gives rise to a colimit-preserving functor $[-]^{t} : \square^{op}\set\to \dtopM$ defined by 
\[
[K]^{t} = \liminj_{\square[n]\rightarrow K} \moore^{top}(|\square[n]|^{t})
\]
\ed

To aid the reader, we recall the definition of spatial precubical set --a notion central to both \cite{RegularMoore} and this addendum.

\begin{nota} \cite[Notation~A.1]{NaturalRealization}
	Let $n\geq 3$. Let $\mathcal{B}_n$ be the set of precubical sets $A$ such that $A\subset \de\square[n]$ and such that $|A|_{geom} \subset [0,1]^n$ contains a $d$-path of $[0,1]^n$ from $(0,0,\dots,0)$ to $(1,1,\dots,1)$ which does not intersect $\{0,1\}^n\backslash\{(0,0,\dots,0),(1,1,\dots,1)\}$. One has $\de\square[n]\in \mathcal{B}_n$.
\end{nota}

\bd \label{def:carac_spatial}
A precubical set is \emph{spatial} if it is orthogonal to the set of maps of precubical sets 
\[\bigg\{\square[n]\sqcup_A \square[n]\longrightarrow \square[n]\mid n\geq 3 \hbox{ and }A\in \mathcal{B}_n\bigg\}.\]
\ed

Every proper precubical set in the sense of \cite[page~499]{MR3722069} is spatial by \cite[Proposition~7.5]{NaturalRealization}. In particular, for all $n\geq 0$, the precubical sets $\de\square[n]$ and $\square[n]$ are spatial, as well as all geometric precubical sets in the sense of \cite[Definition~1.18]{zbMATH07226006} and all non-positively curved precubical sets in the sense of \cite[Definition~1.28]{zbMATH07226006}, since they are proper. Hence many precubical sets arising from real concurrent systems are spatial by \cite[Proposition~1.29]{zbMATH07226006}. Also every $2$-dimensional precubical set is spatial by \cite[Corollary~A.3]{NaturalRealization}.

\bth \label{thm:iso_reg_reg0} 
For all precubical sets $K$, there is a natural weak equivalence of the h-model structure of Moore flows \[[K]^{t} \longrightarrow \moore^{top}(|K|^{t}).\] Moreover, the weak equivalence above is an isomorphism of Moore flows if and only if $K$ is spatial.
\eth

\bpf
It is mutatis mutandis the proof of \cite[Theorem~6.5]{RegularMoore}.
\epf

\bth \label{thm:m-cof-moore-flow}
For all precubical sets $K$, the Moore flow $[K]^{t}$ is m-cofibrant.
\eth

\bpf
It is mutatis mutandis the proof of \cite[Theorem~6.7]{RegularMoore}.
\epf

\begin{cor} \label{cor:space-mcof}
	Let $K$ be a precubical set. Let $(\alpha,\beta)\in K_0\p K_0$. The space of tame $d$-paths from $\alpha$ to $\beta$ in the geometric realization of $K$ is homotopy equivalent to a CW-complex.
\end{cor}

\bpf
Using Theorem~\ref{thm:m-cof-moore-flow}, it is mutatis mutandis the proof of \cite[Corollary~6.8]{RegularMoore}. 
\epf

\bth \label{thm:three} \cite[Theorem~6.14]{QHMmodel} 
Let $r\in \{q,h,m\}$. There exists a unique model structure on $\ptop{\mathcal{M}}$ such that: 
\begin{itemize}
	\item A map of multipointed $d$-spaces $f:X\to Y$ is a weak equivalence if and only if $f^0:X^0\to Y^0$ is a bijection and for all $(\alpha,\beta)\in X^0\p X^0$, the continuous map $\PP_{\alpha,\beta}^{top}X\to \PP^{top}_{f(\alpha),f(\beta)}Y$ is a weak equivalence of the r-model structure of $\Top$.
	\item A map of flows $f:X\to Y$ is a fibration if and only if for all $(\alpha,\beta)\in X^0\p X^0$, the continuous map $\PP^{top}_{\alpha,\beta}X\to \PP^{top}_{f(\alpha),f(\beta)}Y$ is a fibration of the r-model structure of $\Top$.
\end{itemize}
All objects are fibrant. The m-model structure of multipointed $d$-spaces is the mixing of their q-model structure and h-model structure in the sense of \cite[Theorem~2.1]{mixed-cole}. These three model categories are denoted by $\ptop{\mathcal{M}}_q,\ptop{\mathcal{M}}_h,\ptop{\mathcal{M}}_m$.
\eth

\bp 
Let $K$ be a spatial precubical set. Then the multipointed $d$-space $|K|^{t}$ is m-cofibrant if it is h-cofibrant.
\ep

\bpf
Using \cite[Theorem~14]{Moore3} establishing that $\moore^{top}:\ptop{\mathcal{M}}_q\to \dtopM_q$ is a right Quillen equivalence, we deduce that $\moore^{top}:\ptop{\mathcal{M}}_m\to \dtopM_m$ is a right Quillen equivalence as well. Then the proof goes as the one of \cite[Proposition~6.9]{RegularMoore}.
\epf

\end{document}